\documentstyle{article}

\newcommand{\pl}{\partial_}

\newcommand{\ol}{\overline}
\newcommand{\bqq}{\begin{equation} \label}
\newcommand{\eeq}{\end{equation}}
\newcommand{\cc}{{\bf C}}
\newcommand{\rr}{{\bf R}}
\newcommand{\grad}{\,{\rm grad}\,}
\newcommand{\ad}{\,{\rm ad}\,}
\newtheorem{The}{Theorem}

\begin{document}
\mathsurround=2pt

\renewcommand{\thefootnote}{\fnsymbol{footnote}}

\begin{titlepage}
January 5, 1996\hfill{KSU/GRG-96-1}
\vskip3cm
\begin{center}
{\bf \large HAMILTON DYNAMICS AND $H$-PLANAR CURVES}
\vskip2cm

{\large D. A. Kalinin}
\footnote{E-mail: kalinin@phys.ksu.ras.ru}\\
Department of General Relativity \\
and Gravitation, Kazan State University \\
18 Lenin Street, KAZAN, 420008, Russia
\end{center}
\vskip1cm

\begin{abstract}
An important example of Hamilton flows on K\"ahler manifold are
$H$-planar flows: all their trajectories are $H$-planar curves
(complex analog of geodesics). The equation which has to obey the
Hamiltonian of $H$-planar flow is received and the method of
finding general solution of this equation is proposed.

Trajectories of charged particles in magnetic fields of special form
on K\"ahler manifolds of constant holomorphic sectional curvature are
studied. Using the fact that K\"ahler manifolds of constant
holomorphic sectional curvature admit $H$-pro\-jec\-ti\-ve mapping on
flat space $\cc^n$ the equation of particle motion is reduced to an 
ordinary differential equation of second order.
\end{abstract}
\end{titlepage}

\renewcommand{\thefootnote}{\arabic{footnote}}
\setcounter{footnote}{0}
\setcounter{section}{0}


\section{Introduction}

The purpose of this paper is to investigate Hamilton systems on
K\"ahler manifolds from the point of view of differential geometry.
K\"ahler manifolds have exceptional significance in geometry and
physics since they have both Riemannian and symplectic structures and
these two structures are in\-ter\-con\-nec\-ted. Because of this fact
one can consider two dynamics on a K\"ahler manifold $M$: one dynamics
is defined by the Riemannian structure and all its trajectories are
geodesics (or $H$-planar curves \cite{otas,sin1} which can be 
considered as a natural analog of geodesics in complex case) and 
another dynamics is connected with symplectic structure on $M$: it is 
defined by the Hamilton equations.

It is of interest to answer the following question: when this two
dynamics coincides, namely, what are the conditions which Hamiltonian
${\cal H}$ and K\"ahler metric $g$ have to obey in order to 
trajectories of Hamilton flow be geodesics (or in more general 
case $H$- pla\- nar curves). The first part of this paper is devoted 
to the solution of this problem (see Theorem~\ref{th2}, Sec.~2).

At the second part  an particular case of $H$-planar curves which has 
important physical applications are considered. These $H$-planar 
curves can be interpreted as trajectories of charged particles in 
a magnetic field of special form (K\"ahler magnetic fields). Particle 
motion in this case is ruled by the system of $2n$ differential
equation of second order and solution of this system in general case
is a complicated problem. But if considered K\"ahler manifold have
constant holomorphic sectional curvature it can be shown (see Sec.4)
that the number of equations can be reduced to one. It gives the
possibility to investigate the structure of trajectories space by 
qualitative methods or receive the solution numerically.

\medskip
We start from recalling some relevant facts of differential geometry
of K\"ahler and symplectic manifolds \cite{arn1,kobn,sin1}.

A Hermitian manifold $M$ with metric $g$ and complex structure $J$  is
called  {\it  K\"ahler  manifold}  if  its  fundamental  form $\omega$
(defined by the condition $\omega(X,Y)=g(JX,Y)$  for  any  two  vector
fields  $X,Y$  on  $M$) is closed $d\omega =0$.  A K\"ahler manifold 
is said to be  {\it  properly  K\"ahler}  if  its  metric  is  
positively definite  (Riemannian).  Because of integrability of 
complex structure $J$ it is  possible  to  introduce  on  K\"ahler  
manifold  $(M,g,J)$, $\mbox{dim}_{\rr}\   M=2n$   local   complex   
coordinates  $z^\alpha, \ol{z^\alpha}$, $\alpha=1,\ldots,n$ such that 
the components of metric and complex structure have the following 
form \cite{kobn} 
\bqq{compg}
g_{\alpha\ol\beta}=\ol {g_{\ol\alpha\beta}}=
\partial_{\alpha}\partial_{\ol\beta}\Phi, \qquad
g_{\alpha\beta}=g_{\ol\alpha\ol\beta}=0,
\eeq
\bqq{compJ}
J^\alpha_{\beta}=- 
J_{\ol\beta}^{\ol\alpha}=i\delta^\alpha_\beta,\qquad
J^\alpha_{\ol\beta} = J_{\beta}^{\ol\alpha} = 0
\eeq
where bar  denotes complex conjugation and real-valued function $\Phi$
on $M$ is {\it K\"ahler potential} of metric $g$.

Since fundamental  2-form  $\omega$  of  K\"ahler  manifold   $M$   is
non-degenerate  $(M,\omega)$ is {\it symplectic manifold} \cite{arn1}.
For any smooth function $f$ on $M$ its  {\it  gradient  vector  field}
$\grad  (f)$  and  {\it  Hamiltonian  vector field} $\ad (f)$
are defined by the conditions $g(X,  \grad (f)) = df (X)$ and $\omega 
(X,  \ad (f))\omega =df$ where $X$ is a vector field on $M$. It is 
easy to see that vector fields $\grad (f)$ and $\ad (f)$ are 
orthogonal at any point  $p\in M$ and  for  any function $f$:  $g(\ad 
(f),  \grad (f))=-\omega (\ad (f),\grad  (f))=0$.  The  components  
of  these  vector  fields  in  
local coordinates $x^i$ on $M$ can be written in the form
\bqq{grad-ad}
\grad (f)^i =g^{ik} \pl k f, \quad \ad (f)^i = \omega^{ki} \pl k f
\eeq
where $g^{ik}$ and $\omega^{ik}$ are components of matrices inverse to
the matrices $(g_{kj})$, $(\omega_{kj})$ of metric and fundamental
2-form. Dynamical system $\gamma$ defined by {\it Hamilton equation}
\bqq{Ham-eq}
\frac{d\gamma}{dt}=\ad ({\cal H})
\eeq
on $M$ is said to be {\it Hamilton flow} with {\it Hamiltonian} ${\cal
H}$.
Curves $\gamma_t$ are called trajectories of dynamical system
$\gamma$.

A curve $\gamma: [0,1] \to M$ on a K\"ahler $M$ is called {\it
$H$-planar curve} if $\chi\equiv d\gamma / dt$ obeys the equation
\bqq{H-plan}
\nabla_\chi\chi = a(t)\chi +b(t)J(\chi)
\eeq
where $a(t)$, $b(t)$ are some smooth functions of $t$ and $\nabla$ is
Riemannian connection of K\"ahler metric $g$. $H$-pla\-nar curves are
natural analog of geodesics in complex case, namely, $H$-planar curves
are real realiza\-ti\-ons of geodesics in the space over complex
algebra \cite{vishsh}. We shall call Hamilton flow {\it $H$-planar
flow} if all its trajectories are $H$-planar curves.

Let $(M,g,J)$   and   $(M',g',J')$   be   two   K\"ahler    manifolds.
Diffeomorphism   $f:M\to   M'$   is   called   $H$-projective  mapping
\cite{otas}  if  an  image  $f\circ\gamma$  of  any  $H$-planar  curve
$\gamma$  on  $M$  is  $H$-planar  curve on $M'$.  In order to mapping
$f:M\to M'$ be $H$-projective it is necessarily preserving of the
complex structure $J$: $f_* \circ J=J'\circ f_*$ where $f_*$ is
differential of $f$.

For future considerations we shall need the following theorem
\begin{The}\label{th1}
{\rm \cite{sin1}}
K\"ahler manifold $(M,g)$ is $H$-pro\-jec\-ti\-ve\-ly  flat,  i.e.  it
admits  $H$-projective  mapping $f$ on flat space 
$\cc^n\cong\rr^{2n}$,
if and only if its holomorphic sectional curvature $k$ is constant.
\end{The}


\section{$H$-planar flows on a K\"ahler manifolds}

Let us consider K\"ahler manifold $M$ with metric $g$, fundamental
2-form $\omega$ and let $\gamma$ be Hamilton flow on $M$. Since
$H$-planar curves are mostly natural analog of geodesics in complex
case, $H$-planar flows form the important class of Hamilton flows on
K\"ahler manifolds.

The problem arises: {\it what are the conditions which has to obey
Hamiltonian function ${\cal H}$ in order to correspondent Hamilton
flow be $H$-planar} ? The solution of this problem is given by the
following

\begin{The}\label{th2}
Hamilton flow $\gamma$ with Hamiltonian ${\cal H}$ on K\"ahler
manifold $M$ is $H$-planar if and only if ${\cal H}$ obeys the
condition
\bqq{condH}
\nabla_{\grad ({\cal H})} \grad ({\cal H}) = B\grad ({\cal H}) - A
J\grad ({\cal H})
\eeq
where $A$ and $B$ are real-valued functions smooth on $M$. Each
trajectory $\gamma_t$ of the flow $\gamma$ is $H$-planar curve on $M$
satisfying the equation {\rm (\ref{H-plan})} where $a(t) =
A(\gamma_t)$ and $b(t)=B(\gamma_t)$.
\end{The}

\noindent
{\bf Proof.} 
If a Hamilton flow $\gamma$ is $H$-planar then any its trajectory
$\gamma_t$ is $H$-planar curve and obeys the equation (\ref{H-plan}).
At the same time $\gamma$ satisfies the Hamilton equation
(\ref{Ham-eq}). In a local complex coordinates $z^\alpha,
\ol{z^\alpha}$, $\alpha=1,\ldots,n$ on $M$ these two equations take
the form
\bqq{hpl-loc}
\frac{d^2 z^\alpha}{dt^2}+ \Gamma^\alpha_{\mu\nu}\frac{dz^\mu}{dt}
\frac{dz^\nu}{dt} = a(t) \frac{dz^\alpha}{dt} +
b(t) J^\alpha_\mu \frac{dz^\mu}{dt},
\eeq
\bqq{ham-loc}
\frac{dz^\alpha}{dt} = \omega^{\ol\mu\alpha}\pl{\ol\mu} {\cal H}
\eeq
where $\Gamma^\alpha_{\mu\nu}=\ol{\Gamma^{\ol\alpha}_{\ol\mu\ol\nu}}$,
$\Gamma^\alpha_{\ol\mu\nu}=\Gamma^\alpha_{\ol\mu\ol\nu}=
\Gamma^{\ol\alpha}_{\ol\mu\nu}=\Gamma^{\ol\alpha}_{\mu\nu}=0$ are
Cristoffel symbols of metric $g$.

Differentiating (\ref{ham-loc}) we find
$$
\frac{d^2 z^\alpha}{dt^2}=
\pl\nu\omega^{\ol\mu\alpha}
\omega^{\nu\ol\rho} \pl{\ol\mu} {\cal H}\pl{\ol\rho} {\cal H}+
\pl{\ol\nu}\omega^{\ol\mu\alpha}
\omega^{\ol\nu\rho} \pl{\ol\mu} {\cal H}\pl{\rho} {\cal H}+
$$
\bqq{diff}
\omega^{\ol\mu\alpha}\omega^{\nu\ol\rho}\pl{\ol\rho}
{\cal H}\pl{\ol\mu\nu} {\cal H}
+\omega^{\ol\mu\alpha}\omega^{\ol\nu\rho}\pl\rho {\cal H}
\pl{\ol\mu\ol\nu} {\cal H}.
\eeq
 From (\ref{hpl-loc}), (\ref{ham-loc})  and (\ref{diff}) it follows
$$
\Gamma^\alpha_{\mu\nu}\omega^{\mu\ol\rho}\omega^{\nu\ol\tau}
\pl{\ol\rho} {\cal H} \pl{\ol\tau} {\cal H} -
\pl\nu\omega^{\ol\mu\alpha}
\omega^{\nu\ol\rho} \pl{\ol\mu} {\cal H}\pl{\ol\rho} {\cal H} -
$$
$$
\pl{\ol\nu}\omega^{\ol\mu\alpha}
\omega^{\ol\nu\rho} \pl{\ol\mu} {\cal H}\pl{\rho} {\cal H} -
\omega^{\ol\mu\alpha}\omega^{\nu\ol\rho}\pl{\ol\rho}
{\cal H}\pl{\ol\mu\nu} {\cal H}
- \omega^{\ol\mu\alpha}\omega^{\ol\nu\rho}\pl\rho {\cal H}
\pl{\ol\mu\ol\nu} {\cal H}.
$$
\bqq{res}
=a(t) \omega^{\ol\mu\alpha}\pl{\ol\mu} {\cal H} + b(t) J^\alpha_\mu
\omega^{\ol\nu\mu}\pl{\ol\nu} {\cal H} =
(a(t)+ib(t)) \omega^{\ol\mu\alpha}\pl{\ol\mu} {\cal H}.
\eeq
Using the identities $\omega_{\alpha\ol\beta}=i g_{\alpha\ol\beta}$,
$\Gamma^\alpha_{\mu\nu}=g^{\alpha\ol\rho}\pl\mu g_{\nu\ol\rho}=
-\omega^{\alpha\ol\rho}\pl\mu \omega_{\nu\ol\rho}$ and
Eqs.~(\ref{compg}), (\ref{compJ}) we can write (\ref{res}) in the
following form
\bqq{loc-cond}
{\cal H}_{,\alpha\ol\mu} g^{\rho\ol\mu} {\cal H}_{,\rho}+
{\cal H}_{,\alpha\mu} g^{\mu\ol\rho} {\cal H}_{,\ol\rho}=
(b - ia) {\cal H}_{,\alpha}.
\eeq
where comma denotes covariant differentiation with respect to the
metric $g$. Taking into account the reality of the function ${\cal H}$
and Eqs.~(\ref{compg}) -- (\ref{grad-ad}) we obtain (\ref{condH})
where functions $A$ and $B$ coincide with $a(t)$, $b(t)$ for any
trajectory $\gamma_t$ of $H$-planar flow $\gamma$.

Let we have a general solution of the equation
$$
\nabla_\chi \chi = b(t)\chi - a(t) J\chi, \quad
\chi=\frac{d\gamma}{dt}
$$
of $H$-planar curve in a $2n$-dimensional K\"ahler manifold $M$:
\bqq{hpl-sol}
x^i=x^i(t,c^l_1,c^l_2), \quad i,l=1,\ldots,2n
\eeq
where $c^l_1,c^l_2$ are integration constants.

Now we propose the method of finding general solution of 
Eq.~(\ref{condH}).

If we take $2n+1$ numbers of $c^l_1,c^l_2$, $l=1,\ldots,2n$ to be
constant and vary rest $2n-1$ numbers we obtain a family of $H$-planar
curves in a open domain $U\subset M$
\bqq{fam}
x^i=x^i(t,\sigma^A), \quad A=1,\ldots,2n-1.
\eeq
Let us calculate the time derivative of Hamiltonian ${\cal H}$ along
the curves of this family
$$
\frac{d{\cal H}}{dt}=\frac{\partial {\cal H}}{\partial x^i}
\frac{dx^i}{dt}=g_{ik}\frac{dx^i}{dt}\frac{dx^k}{dt}
$$
then
$$
{\cal H}=\int dt \ g_{ik}\frac{dx^i}{dt}\frac{dx^k}{dt} + {\rm const}.
$$
By using (\ref{fam}) we can eliminate $t$ and $\sigma^A$,
$A=1,\ldots,2n-1$ from these expression and obtain Hamiltonian ${\cal
H}$ as a function of coordinates $x^i$, $i=1,\ldots,2n$. Considering
different families of $H$-planar curves of the form (\ref{fam}) we
receive different Hamiltonians obeying (\ref{condH}).

Since any Hamilton flow with a Hamiltonian ${\cal H}$ defines
a $2n-1$-pa\-ra\-me\-ter family of $H$-planar curves we described the
way of construction general solution of Eq.~(\ref{condH}) using
general solution of equation of $H$-planar curve in a $2n$-dimensional
K\"ahler manifold $M$.

In the conclusion of this section let us consider gradient
flow $\gamma_{gr}$ of the function ${\cal H}$:
\bqq{11*}
\frac{d\gamma}{dt}=\grad{{\cal H}}.
\eeq
Writing (\ref{loc-cond}) for trajectories of the system (\ref{11*}) we
find
$$
\nabla_\chi \chi = b(t)\chi + a(t) J\chi, \quad
\chi=\frac{d\gamma}{dt}.
$$
 From here it is seen that considering trajectories are $H$-planar
curves. It yields the following
\begin{The}\label{th2-1}
Let $\gamma$ be $H$-planar Hamilton flow with Hamiltonian ${\cal H}$ 
on
K\"ahler manifold $M$. Then trajectories of gradient flow
$\gamma_{gr}$ of function ${\cal H}$ are $H$-planar curves.
\end{The}


\section{Magnetic fields on Riemannian manifolds}

Magnetic flows on (pseudo) Riemannian manifolds are interesting
example of dynamical systems whose trajectories (trajectories of
charged particles in magnetic field) in general case are not geodesics
but closely related with Riemannian structure of the manifold. If
considered manifold is K\"ahler one there is an interesting example of
magnetic fields for which these trajectories are $H$-planar curves ---
K\"ahler magnetic fields.

Let $(M,g)$ be $m$-dimensional pseudo Riemannian manifold. Closed
2-form $B$ on $M$ is said to be {\it magnetic field}. If $M$ is
Lorentzian manifold then magnetic fields on $M$ usually refers to as
{\it electromagnetic}.

Let ${\cal  I}$  be  smooth field of skew-symmetric linear operators
${\cal I}_p: T_p M \to T_p M$ defined by the condition 
\bqq{5s1-1}
B(X, Y) = g ({\cal I}X, Y)
\eeq
for any vector fields $X, Y$ on $M$.
A smooth curve $\gamma: I=[0,1]\to M$ is called {\it trajectory of
magnetic field} $B$ if its tangent vector $\chi \equiv d\gamma (t) /
dt$ satisfies the equation
\bqq{5s1-2}
\nabla_ {\chi} \chi= {\cal I} (\chi).
\eeq

Since operator  ${\cal I}$ is skew-symmetric, any trajectory $\gamma$
of magnetic field $B$ on $M$ has constant speed, i.e.
\bqq{5s1-speed}
\frac{d}{dt} g(\chi,  \chi) = 2g({\cal I}(\chi), \chi) = 0.
\eeq
If $\gamma (t)$ is trajectory of magnetic field $B$ with the speed $V$
then $\gamma ( \alpha t) $, $\alpha \in \rr$ is trajectory of magnetic
field $\alpha B$ with the same speed.

The magnetic  field  $B$  on pseudo Riemannian manifold $( M,  g) $ is
called {\it uniform} \cite{ad1,ad2,comtet} if tensor field ${\cal I} 
$  is
parallel with  respect to Riemannian connection $\nabla$ of the metric
$g$:  $\nabla {\cal I} =0$.  In the case of  flat  manifold  $M$  such
definition  coincides  with  the  definition of a homogeneous magnetic
field in the flat space.

Let $(M,g,J)$ be a K\"ahler manifold of real dimension $2n$ and
$\omega$ is fundamental 2-form of this manifold. Closed 2-form
$q\omega$ where $q\neq 0$ is a real constant is uniform magnetic field
which we call {\it K\"ahler magnetic field} on $M$ \cite{ad1,ad2}.

Trajectories of uniform magnetic fields on surfaces of constant
curvature have studied Comtet \cite{comtet} and Sunada \cite{sunada}.
They have proven that on the sphere $S^2$ and on Euclidean plane
$\rr^2$ these trajectories are circles (in the first case they are
circles in the sense of Euclidean geometry of $\rr^3$ in which $S^2$
is canonically imbedded). On the hyperbolic plane of sectional
curvature $-k$ the behaviour of trajectories principally differs. If
$|q|\ge\sqrt{k}$, then unit speed trajectories are still closed curve,
but when $|q|\le\sqrt{k}$ the trajectories of unit speed are open
curve. Adachi \cite{ad1,ad2} has generalized these results on the 
case of properly K\"ahler manifolds of any dimensions.

The studies of K\"ahler magnetic fields on pseudo Ri\-e\-man\-ni\-an
manifolds is of interest from the point of view of applications in
theoretical physics, because curves described by the equation
(\ref{5s1-2}) in four-dimensional manifold can be interpreted as
trajectories of charged particles in electromagnetic and
gra\-vi\-ta\-ti\-o\-nal fields. Similar curves in pseudo Riemannian
and, in particular, K\"ahler manifolds of any signature and dimension
arise as trajectories of particles in various field theories
\cite{vl-yu-s,ing,pav}.

In the next section we investigate trajectories of K\"ahler magnetic
fields on a K\"ahler manifold of constant holomorphic sectional
curvature and any signature and dimension (see also \cite{ka1,ka2}).


\section{K\"ahler magnetic fields
\newline
on manifolds of constant
\newline
holomorphic sectional curvature}

Let $M$  be  K\"ahler  manifold  of  constant  holomorphic   sectional
curvature $k$ with K\"ahler metric $g$ and fundamental 2-form $\omega$
and let $B=q\omega$, $q\in \rr$ is K\"ahler magnetic field on $M$. The
trajectory  $\gamma_t$ of magnetic field $B$ satisfies the equation
(\ref{5s1-2})     which     in     local     complex      coordinates
$z^\alpha$,$\ol{z^\alpha}$,  $\alpha=1,\ldots,n$ on $M$ takes the form
\bqq{5s1-3}
\nabla_t \frac {d z^\alpha}{dt}  =  iq  \frac{dz^\alpha}{dt},  \qquad
\nabla_t \frac{d\ol{z^\alpha}}{dt}
=-iq\frac{d\ol{z^\alpha}}{dt}
\eeq
where $\nabla$ is covariant derivative with respect to the metric  $g$
and $\nabla_t \equiv dz^\sigma / dt \,\nabla_\sigma +  d\ol{z^\sigma}/
dt\, \nabla_{\ol\sigma}$.

According to Theorem~\ref{th1} it  exists
$H$-pro\-jec\-ti\-ve     mapping     $f:\cc^n     \to     M$.    Since
$H$-pro\-jec\-ti\-ve mapping preserve  the  complex  structure  it  is
possible  to  choose  on  $\cc^n$  complex  coordinates  in  which all
Cristoffel symbols of the metric are equal to zero.  Then in  
corresponding
local complex coordinates on $M$ Cristoffel symbols are defined
by the following formula
\bqq{5s1-3b}
\Gamma^\lambda_{\alpha\mu}=
\psi_\alpha \delta^\lambda_\mu + \psi_\mu \delta^\lambda_\alpha
\eeq
where $\psi_{\alpha}=\pl\alpha\psi$ and $\psi$ is a smooth function on
$M$. K\"ahler manifolds of constant holomorphic sectional curvature
are Einstein manifolds and, hence \cite{kobn},
\bqq{5s1-3bb}
R_{\alpha\ol\beta} =\frac{1}{2}k(n+1) g_{\alpha\ol\beta},
\quad R_{\alpha\beta}=0, \quad R_{\ol\alpha\ol\beta} =0,
\eeq
where $R_{ij}$ are components of Ricci tensor. From the equalities
$R_{\alpha\ol\beta} =
\partial_{\ol\beta}\Gamma^\sigma_{\alpha\sigma}$,
$\ol{R_{\alpha\ol\beta}} = R_{\ol\alpha\beta}$ (see, for example,
\cite{kobn}) using (\ref{5s1-3b}) and (\ref{5s1-3bb}) we obtain
$$
\partial_{\alpha\ol\beta}\psi
=-\frac{k}{4}\partial_{\alpha\ol\beta} \Phi.
$$
Hence,
$$
\psi=-\frac{k}{4} \Phi
$$
up to the function of the form $\varphi (z) + \ol \varphi (\ol z)$
which is omitted because it corresponds to admissible transformations 
of
K\"ahler of potential $\Phi\to\Phi+\frac{k}{4}(\varphi (z)+\ol{\varphi
(z)})$.
Using Theorem~\ref{th1} from here  we  find
\bqq{5s1-zzz}
\psi= -\frac{1}{2}\ln(\sum\limits_\alpha
\epsilon_\alpha  z^\alpha  \ol{z^{\alpha}}) + 1, \quad \epsilon_\alpha
=\pm 1.
\eeq

Let $\gamma$ be $H$-planar curve determined by Eq.~(\ref{5s1-3}) and
$f:M\to M'$ be $H$-projective mapping. In respective complex
coordinates the curve $f^{-1}\gamma$ is defined by the same equation
as $\gamma$ (i.e. $z^\alpha = z^\alpha (t) $, $\ol{z^\alpha}
=\ol{z^\alpha}(t)$, $\alpha=1,\ldots,n$) and satisfies the equation of
$H$-planar curve in $\cc^n$:
\bqq{5s1-4}
\frac{d^2 z^\alpha}{dt^2}=c(t)\frac{dz^\alpha}{dt},
\qquad\frac{d^2\ol{z^\alpha}}{dt^2}=\ol{c(t)}\frac{d\ol{z^\alpha}}{dt}
\eeq
where $c(t)$ is a smooth complex-valued function.

Thus, the dynamical system (\ref{5s1-3}) can be simulated dynamical
system (\ref{5s1-4}). Solving equation (\ref{5s1-4}) we find.
\bqq{5s1-5}
z^\alpha= C_1^\alpha f (t) + C^\alpha_2
\eeq
where $f (t) =\int dt(\exp\int c(t) dt)$, and $C_1^\alpha$,
$C^\alpha_2$ are integration constants.

Let us consider two cases: $\sum\limits_\alpha\epsilon_\alpha
C_1^\alpha \ol{C_1^\alpha} \neq 0$ and
$\sum\limits_\alpha\epsilon_\alpha C_1^\alpha \ol { C_1^\alpha } = 0$.
In the first case making in (\ref{5s1-5}) the following transformation
\bqq{5s1-6}
\tilde f (t) = vf (t) + w, \quad v,w\in \cc
\eeq
and presenting  $z^\alpha$  in 
 the form $ z^\alpha= \tilde C_1^\alpha
\tilde f (t) + \tilde C^\alpha_2$,  it is possible to choose  
constants
$v$, $w$ in such a way that
\bqq{5s1-7}
\sum\limits_\alpha\epsilon_\alpha \tilde C_1^\alpha \ol { \tilde
C_2^\alpha   }  =  0,\qquad  \sum\limits_\alpha\epsilon_\alpha  \tilde
C_1^\alpha \ol{\tilde C_1^\alpha} = {\cal A},
\eeq
where ${\cal A}= {\rm sgn}\, (\sum\limits_\alpha\epsilon_\alpha
C_1^\alpha \ol{C_1^\alpha})$.

By the same way in the case $\sum_\alpha\epsilon_\alpha C_1^\alpha
\ol{C_1^\alpha} = 0$ making transformation (\ref{5s1-6}) one can
choose constant $v$ and $w$ so that
\bqq{5s1-8}
\sum\limits_\alpha\epsilon_\alpha \tilde C_2^\alpha \ol { \tilde
C_2^\alpha  }  =-1,  \qquad  \sum\limits_\alpha\epsilon_\alpha  \tilde
C_1^\alpha \ol { \tilde C_2^\alpha } =1.
\eeq

Hereinafter we shall consider that conditions (\ref{5s1-7}) or
(\ref{5s1-8}) holds and tilde over $C^\alpha_1$, $C^\alpha_2$ and $f$
will be omitted.

Let us present complex valued function $f(t)$ in the following form
\bqq{5s1-9}
f(t) = \exp (r(t) + i\phi (t)),
\eeq
where $r(t)=\ol{r(t)}$ and $\phi (t) = \ol { \phi (t) } $. From
(\ref{5s1-6}) it follows that trajectories $\gamma_t$ of K\"ahler
magnetic field $B=q\omega$ belong to 2-dimensional submanifold
$P\subset M$ defined in the complex coordinates $z^\alpha$,
$\alpha=1,\ldots, n$ by the equations
\bqq{5s1-10}
z^\alpha=C_1^\alpha u + C_2^\alpha,\qquad u\in \cc.
\eeq
Thus, the functions $r$ and $\phi$ can be treated as coordinates of
the trajectory's point on the surface $P$. Since any geodesics is a
trajectory of K\"ahler magnetic field $B_0=q_0 \omega$, $q=0$ it is
easy to see that $P$ is properly geodesic submanifold in $M$. Using
the formula \cite{kobn}
$$
\Gamma^\alpha_{\beta\gamma}= \ol{\Gamma^{\ol\alpha}_{\ol\beta 
\ol\gamma}}
=g^{\alpha\ol\mu}\partial_\beta   g_{\ol\mu\gamma},
$$
from (\ref{5s1-3}) we receive
\bqq{5s1-11}
\ddot z^\alpha= - 4\partial_\nu\psi\dot z^\nu\dot z^\alpha + iq\dot
z^\alpha
\eeq
where dot denotes derivation on parameter $t$.

So we obtain $n$ complex equations which can be rewritten as $2n$ real
ones
$$
\frac{d}{dt}\ln(\dot z^\alpha\dot z^{\ol\alpha})=-4\frac{d\psi}{dt},
$$
\bqq{5s1-12}
i \frac{d}{dt}\ln\frac{\dot z^\alpha}{\dot z^{\ol\alpha}} =
-4i(\partial_\nu\psi \dot z^\nu- \partial_{\ol\nu}\psi\dot z^{\ol\nu})
-2 q,
\eeq
where the summation on $\alpha$ is not carried out.

In the   case   $\sum\limits_\alpha\epsilon_\alpha  \tilde  C_1^\alpha
\ol{\tilde C_2^\alpha}=0$,  $\sum\limits_\alpha\epsilon_\alpha  \tilde
C_1^\alpha \ol{\tilde C_1^\alpha}= {\cal A}$ from here with the
help of (\ref{5s1-zzz}), (\ref{5s1-5}) and (\ref{5s1-9}) we receive
\bqq{5s1-13}
\frac{d}{dt}\ln(e^{2r}(\dot r^2 +\dot\phi^2))=-4
\frac{d\psi}{dt},
\eeq
\bqq{5s1-14}
i\frac{d}{dt}\ln
\frac{\dot r + i\dot \phi}{\dot  r-i\dot\phi}-2\dot\phi=
-4i(\partial_\nu\psi \dot z^\nu-\partial_{\ol\nu}\psi
\dot z^{\ol\nu})-2q,
\eeq

Using Eqs.~(\ref{5s1-zzz}),    (\ref{5s1-5}),    (\ref{5s1-7})     and
(\ref{5s1-9}) from (\ref{5s1-13}), (\ref{5s1-14}) we find
\bqq{5s1-15}
\frac{d}{dt}\ln (e^{2r}(\dot r^2+ \dot \phi^2))=
2\frac{d}{dt}\ln ({\cal A}e^{2r}+  C+1),
\eeq
\bqq{5s1-16}
i \frac{d}{dt}\ln \frac{\dot r + i\dot \phi}{\dot r - i\dot \phi}-
2\dot \phi=-4({\cal A}e^{2r}+C+1)^{-1} {\cal A} e^{2r}\dot\phi-2q,
\eeq
where $C=\sum\limits_\alpha\epsilon_\alpha C_2^\alpha
\ol{C_2^\alpha}$.
Integrating (\ref{5s1-15}) we receive
\bqq{5s1-17}
{\cal J}^2 e^{2r}(\dot r^2+\dot\phi^2)=({\cal A}e^{2r}+C+1)^2
\eeq
where ${\cal J}$ is an integration constant.

Using (\ref{5s1-4}) it is possible to calculate a speed square  $V=  2
g_{\alpha\ol\beta}  \dot  z^\alpha  \dot  z^ { \ol\beta } $ of
trajectory $z^\alpha (t) $
$$
V= \frac{4}{k}\frac{(\dot r^2+ \dot \phi^2)e^{2r}(C+1)}{({\cal A}
e^{2r}+ C+1)^2}.
$$
 From here with the help of (\ref{5s1-17}) it follows
$$
V= \frac{4}{k} \frac{C + 1}{{\cal J}^2},
$$
where $k$ is holomorphic sectional curvature of manifold $M$.
According to (\ref{5s1-speed}) any trajectory $\gamma$ of magnetic
field $B$ has constant speed, hence, $V=$const. So we find that {\em
for a trajectory of K\"ahler magnetic field with fixed value of $C$
the integration constant ${\cal J}$ is inverse proportional to square
root of speed}.

Let us consider the coordinate $r$ as function of coordinate $\phi$ on
the surface $P\in M$: $r=r(\phi)$. We choose $\phi$ to be an
independent variable and denote
$$
p\equiv \frac{d\phi}{dt},\quad  r'\equiv\frac{dr}{d\phi}
\quad {\rm then} \quad \frac{dr}{dt}= 
\frac{dr}{d\phi}\frac{d\phi}{dt}=
r'p.
$$
Then Eqs.~(\ref{5s1-15}) and (\ref{5s1-16}) takes the following form
\bqq{5s1-18}
2p(\frac{r''}{{r'}^2+1}- 1)= -4({\cal A}e^{2r}+ C+1)^{-1}{\cal A}
e^{2r} p -2q,
\eeq
$$
e^{2r} p^2 ({r'}^2+1){\cal J}^2= ({\cal A}e^{2r}+ C+1)^2.
$$
Expressing $p$ from the last equation and substituting it in 
(\ref{5s1-18})
we receive
$$
\frac{r''}{r'^2+1} = 1-(Ae^{2r}+C+1)^{-1} 2Ae^{2r}
$$
\bqq{5s1-20}
- q(Ae^{2r}+C +1)^{-1}{\cal J}e^r(r'^2+1)^{1/2}.
\eeq

Similarly, in the case when $\sum_\alpha\epsilon_\alpha
C_1^\alpha\ol{C_1^\alpha} = 0$ from (\ref{5s1-zzz}), (\ref{5s1-5}),
(\ref{5s1-9}) (\ref{5s1-13}) and (\ref{5s1-14}) we find
\bqq{5s1-21}
\frac{r''}{r'^2+1} = -r'\,{\rm tg}\,\phi-
q{\cal J}(1+r'^2)^{1/2}(\cos\phi)^{-1}.
\eeq

\begin{The}\label{th3}
Let $(M,g,J)$ be geodesically complete K\"ahler manifold of constant
holomorphic sectional curvature $k$ with fundamental {\rm 2}-form
$\omega$ and $B=q\omega$, $q\in \rr$ be K\"ahler magnetic field on
$M$. Let $z^\alpha, \ol{z^\alpha}$, $\alpha=1,\ldots,n$ be such local
complex coordinates on $M$ that Cristoffel symbols of metric $g$ have
the form {\rm (\ref{5s1-3b})}. Then any trajectory $\gamma (t)$ of
magnetic field $B$ belongs to {\rm 2}-dimensional properly geodesic
submanifold $P\subset M$ determined in coordinates $(z^\alpha,
\ol{z^\alpha})$ by the equations
$$
z^\alpha=C_1^\alpha u+
C_2^\alpha,\qquad u\in \cc.
$$

The coordinates $r,\phi$ of the point $p=\gamma_t \in P  \subset  M$
satisfy  the  equation  {\rm(\ref{5s1-20})} if ${\cal A}= 
\sum\limits_\sigma
C_1^\sigma \ol{C_1^\sigma} \neq 0$ and equation {\rm(\ref{5s1-21})} 
if $\sum\limits_\sigma C_1^\sigma\ol{C_1^\sigma} = 0$.
\end{The}

This theorem reduces the study of of equations (\ref{5s1-3}) to the
consideration of {\em one ordinary differential equation} of the 
second order. The solution of this equation can be obtained by 
numerical methods.

\section*{Acknowledgments}

Author would like to thank Prof. A.V.~Aminova for stimulating
discussions. The work was supported by Grant No~a496-m of Soros
Science Education Program.


\end{document}